\documentclass[11pt, twoside]{article}

\usepackage[text={390pt, 592pt}, centering, headheight=13.75pt]{geometry}
\usepackage{fancyhdr}
\usepackage{amssymb, amsmath}
\usepackage{graphicx}
\usepackage{hyperref}
\usepackage{breakurl}
\usepackage[labelsep=none]{caption}

\renewcommand{\headrulewidth}{0.5pt}

\hypersetup{pdfstartview=FitH, pdfpagemode=UseNone}

\newcommand{\floor}[1]{\left \lfloor #1 \right \rfloor}
\newcommand{\ceil}[1]{\left \lceil #1 \right \rceil}

\allowdisplaybreaks

\begin{document}

\vspace*{\bigskipamount}

\noindent{\LARGE \bfseries A Proof of Willcocks's Conjecture}\bigskip

\noindent{\scshape Nikolai Beluhov}\bigskip

\begin{center} \parbox{352pt}{\footnotesize \emph{Abstract}. We give a proof of Willcocks's Conjecture, stating that if $p - q$ and $p + q$ are relatively prime, then there exists a Hamiltonian tour of a $(p, q)$-leaper on a square chessboard of side $2(p + q)$. The conjecture was formulated by T.\ H.\ Willcocks in 1976 and has been an open problem since.} \end{center}

\section{Introduction}

A \emph{$(p, q)$-leaper} is a fairy chess piece generalising the knight. On a rectangular chessboard, it can leap from a cell $(x, y)$ to any of the cells $(x \pm p, y \pm q)$ and $(x \pm q, y \pm p)$.

The construction of Hamiltonian tours of the knight has fascinated puzzle-solvers since at least the ninth century (\cite{J1}).

\begin{figure}[!h] \centering \includegraphics{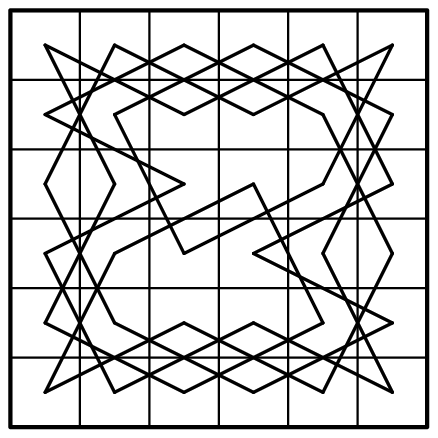} \caption{} \label{euler} \end{figure}

The mathematical study of knight tours began with Leonhard Euler's paper \cite{E}, one of the earliest papers on combinatorics. In particular, Euler gave the first Hamiltonian tour of the knight on a board of size $6 \times 6$ (Figure \ref{euler}), the smallest square board admitting such a tour. Euler's tour possesses central symmetry.

\begin{figure}[!p] \centering \includegraphics{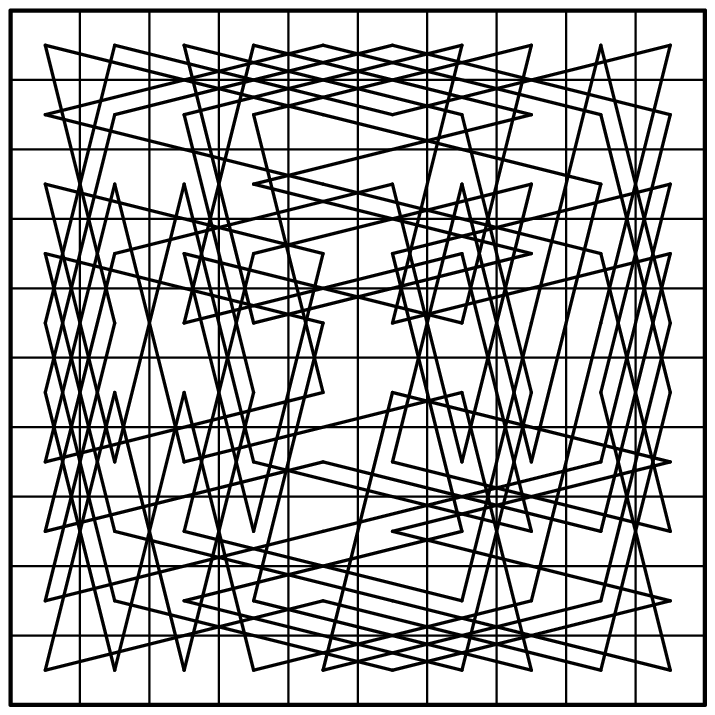} \caption{} \label{frost-1} \end{figure}

\begin{figure}[!p] \centering \includegraphics{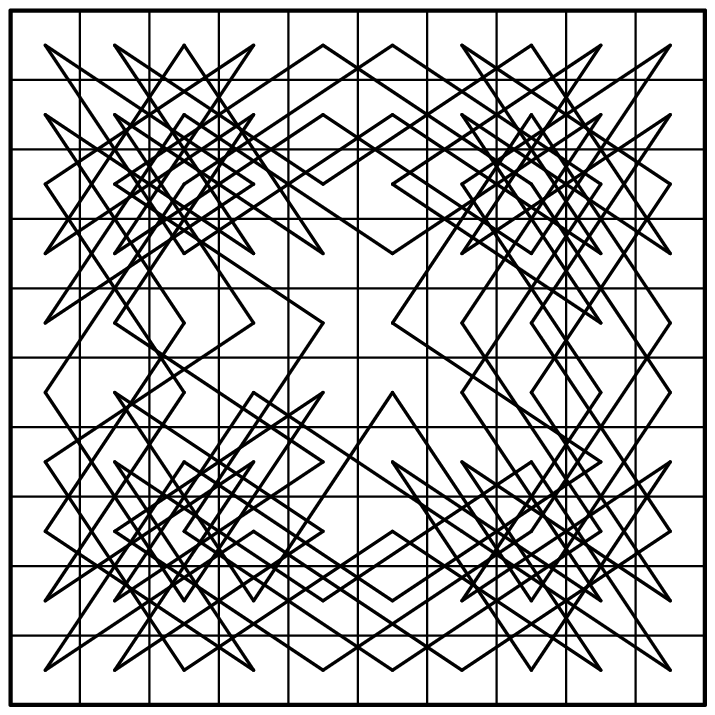} \caption{} \label{frost-2} \end{figure}

The generalisation to leapers followed in the nineteenth century. A.\ H.\ Frost was the first to give Hamiltonian tours of the $(1, 4)$-leaper (the giraffe) and the $(2, 3)$-leaper (the zebra) on a board of size $10 \times 10$ (Figures \ref{frost-1} and \ref{frost-2}), in \cite{F}.

\begin{figure}[!p] \centering \includegraphics{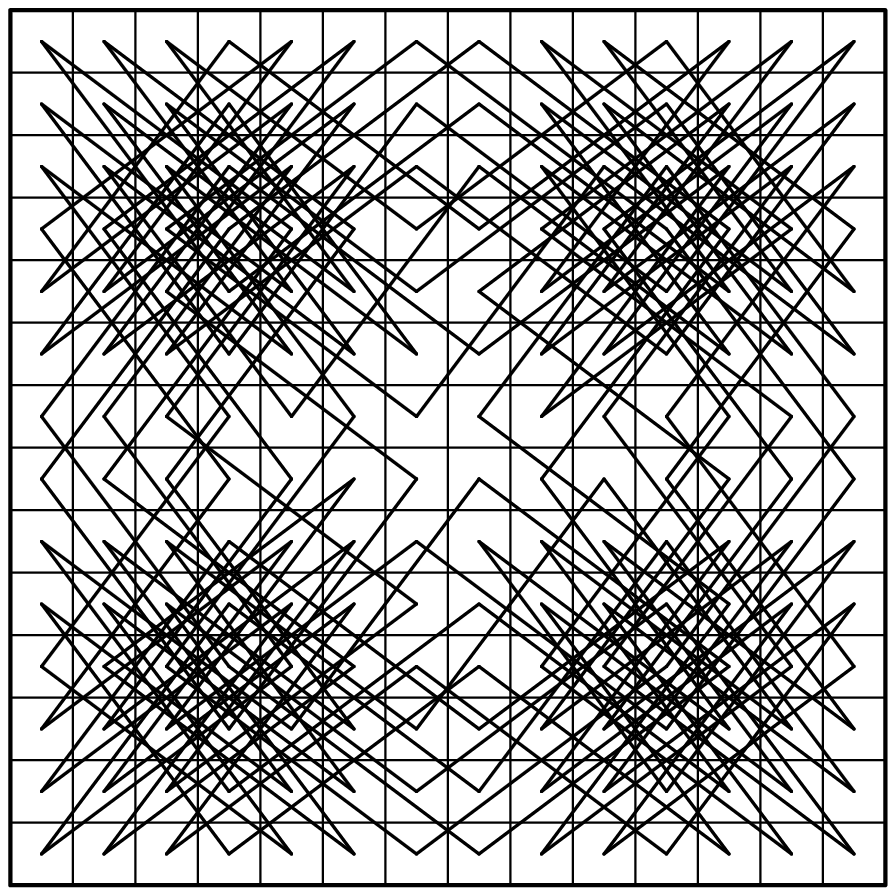} \caption{} \label{willcocks-1} \end{figure}

\begin{figure}[!p] \centering \includegraphics{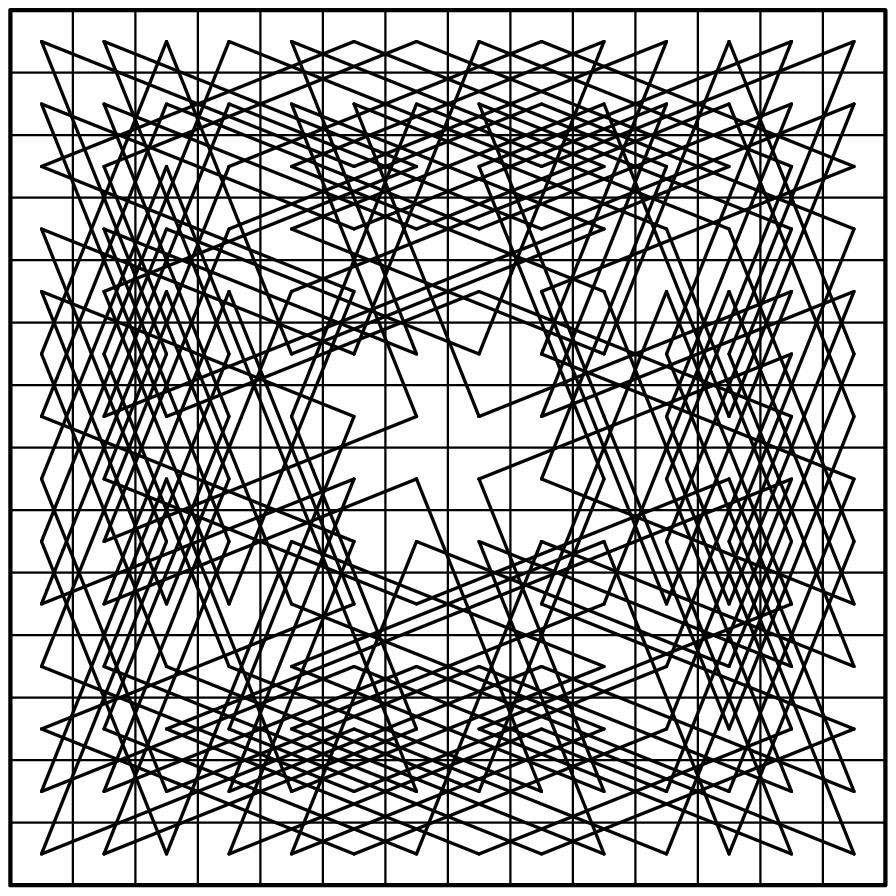} \caption{} \label{willcocks-2} \end{figure}

Nearly one century later, Theophilus H.\ Willcocks gave the first Hamiltonian tours of the $(3, 4)$-leaper (the antelope) and the $(2, 5)$-leaper on a board of size $14 \times 14$ (Figures \ref{willcocks-1} and \ref{willcocks-2}), in \cite{JW}. Both of Willcocks's tours possess fourfold rotational symmetry.

In all cases, those are the smallest square boards admitting such tours, with or without symmetry.

A necessary condition for a Hamiltonian tour of a $(p, q)$-leaper $L$ to exist on a chessboard $B$ is that $L$ can reach from any cell of $B$ to any other by a series of moves. In that case, $L$ is \emph{free on $B$}. A leaper that is free on a rectangular board of size greater than $1 \times 1$ is also free on all larger boards. Such leapers are known as \emph{free}.

The notion of a free leaper was introduced by George P.\ Jelliss and Willcocks in \cite{JW}, where they show that a necessary condition for a leaper to be free is that $p - q$ and $p + q$ are relatively prime. Jelliss conjectured in \cite{J2}, and Donald E.\ Knuth proved in \cite{K}, that this condition is also sufficient.

Willcocks conjectured in \cite{JW} that if a $(p, q)$-leaper $L$ is free, then there exists a Hamiltonian tour of $L$ on a square board of side $2(p + q)$. By the time that Willcocks formulated his conjecture, its truth had been established for the knight, giraffe, zebra, antelope, and $(2, 5)$-leaper.

In \cite{K}, Knuth proved a related conjecture of Willcocks (\cite{JW}), that no square board of side less than $2(p + q)$ admits a Hamiltonian tour of $L$, and verified the existence of a Hamiltonian tour of $L$ on a board of side $2(p + q)$ in the special cases $p = 1$ and $q = p + 1$.

We proceed to give a proof of Willcocks's Conjecture.

\section{Preliminaries}

First we lay some groundwork.

The \emph{leaper graph} of a leaper $L$ on a board $B$ is the graph whose vertices are the cells of $B$ and whose edges join the cells of $B$ joined by a move of $L$. (\cite{J2, K})

In all that follows, let $p$ and $q$ be positive integers such that $p < q$ and $p - q$ and $p + q$ are relatively prime, let $L$ be a $(p, q)$-leaper, and let $B$ be a square board of side $2(p + q)$.

Let $W$, the \emph{Willcocks graph} of $L$, be the leaper graph of $L$ on $B$.

Introduce a Cartesian coordinate system such that the lower left corner of $B$ is at the origin and the upper right corner of $B$ is at $(2(p + q), 2(p + q))$. We refer to each cell of $B$ by the coordinates of its lower left corner, so that the leftmost lowermost cell of $B$ is $(0, 0)$ and the rightmost uppermost one is $(2(p + q) - 1, 2(p + q) - 1)$.

We write $[x_1, x_2] \times [y_1, y_2]$ for the subboard of $B$ whose lower left corner is at the point $(x_1, y_1)$ and whose upper right corner is at the point $(x_2, y_2)$. It consists of all cells $(x, y)$ such that $x_1 \le x < x_2$ and $y_1 \le y < y_2$.

We say that the cell $(x_1 + x, y_1 + y)$ is at \emph{position $(x, y)$} in the subboard $[x_1, x_2] \times [y_1, y_2]$ of $B$.

All of our constructions are going to be highly symmetric. For this reason, we introduce shorthand for the copies of objects within $B$ under the symmetries of $B$.

Let $A$ be any object within $B$ consisting of cells, such as a subboard of $B$ or a path in $W$. Then the \emph{reflections} of $A$ are $A$ and the images of $A$ under reflection in the axis of symmetry $x = p + q$ of $B$, the center of $B$, and the axis of symmetry $y = p + q$ of $B$.

In general, an object $A$ has four reflections. However, some of them may coincide if $A$ is itself symmetric.

Given a cell $a = (x, y)$ and an integer vector $v = (x_v, y_v)$, we write $a + v$ for the cell $(x + x_v, y + y_v)$.

We refer to the eight vectors $(\pm p, \pm q)$ and $(\pm q, \pm p)$ as the \emph{directions} of $L$. The direction of a move of $L$ from $a$ to $b$ is the unique direction $d$ such that $a + d = b$.

Given a cell $a$ and a direction $d$, we write $a \to d$ for the edge of $W$ joining $a$ and $a + d$. Given a sequence of directions $d_1$, $d_2$, $\ldots$, $d_k$, we write $a \to d_1 \to d_2 \to \ldots \to d_k$ for the path in $W$ starting from $a$ whose successive moves are of directions $d_1$, $d_2$, $\ldots$, $d_k$. The vertices of this path are, in this order, $a$, $a + d_1$, $a + d_1 + d_2$, $\ldots$, $a + d_1 + d_2 + \ldots + d_k$.

Similarly, given a subboard $S$ of $B$ and a direction $d$, we write $S \to d$ for the \emph{pencil} of all edges of $W$ of the form $a \to d$ for $a$ in $S$. Given a sequence of directions $d_1$, $d_2$, $\ldots$, $d_k$, we write $S \to d_1 \to d_2 \to \ldots \to d_k$ for the pencil of all paths in $W$ of the form $a \to d_1 \to d_2 \to \ldots \to d_k$ for $a$ in $S$.

\emph{Halving} an even-length cycle $C$ consists in deleting every second edge of $C$.

A \emph{two-factor} of a graph $G$ is a spanning subgraph of $G$ in which every vertex is of degree two. A two-factor is always the disjoint union of several cycles. A connected two-factor is a Hamiltonian tour.

In the context of leaper graphs, a two-factor is also known as a \emph{pseudotour} (\cite{J1}). We are going to use general graph-theoretic and leaper-specific vocabulary interchangeably.

\section{A Proof of Willcocks's Conjecture}

Our proof proceeds in two steps. First we show that the Willcocks graph admits a pseudotour by constructing a family of such pseudotours. Then we show that there exists a member of that family which is in fact a Hamiltonian tour. The second step is more difficult than the first one, and is itself broken into several substeps.

\begin{figure}[!p] \centering \includegraphics{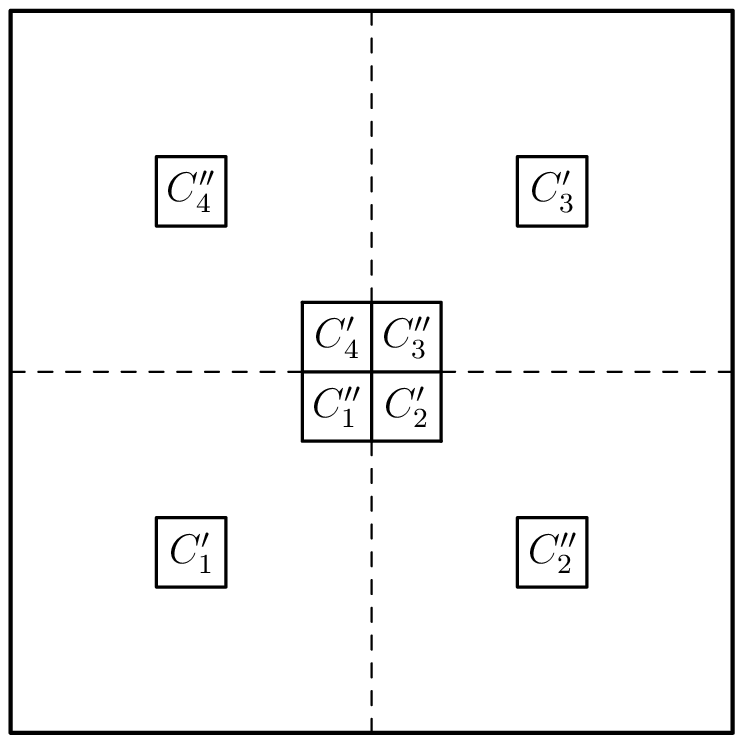} \caption{} \label{i-1} \end{figure}

\begin{figure}[!p] \centering \includegraphics{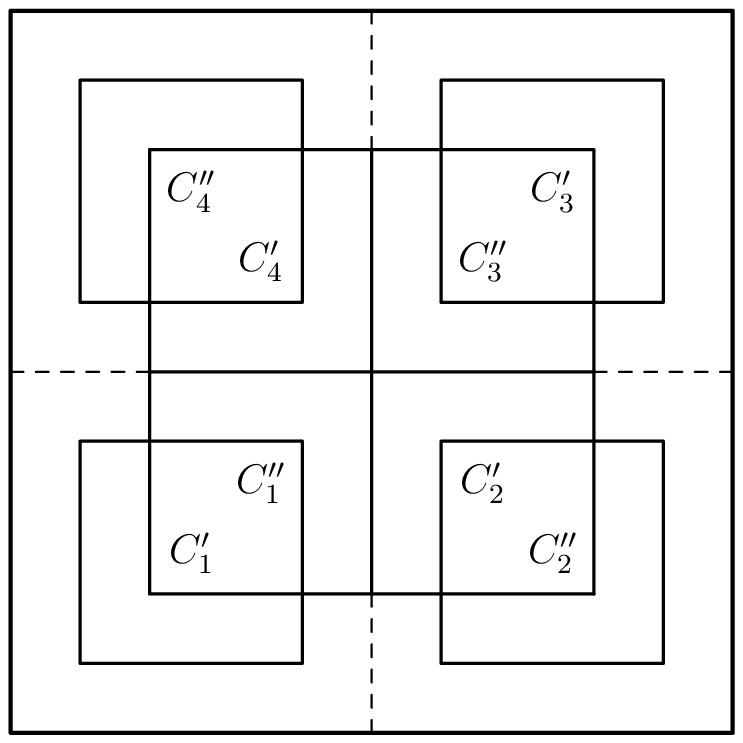} \caption{} \label{i-2} \end{figure}

\begin{figure}[!p] \centering \includegraphics{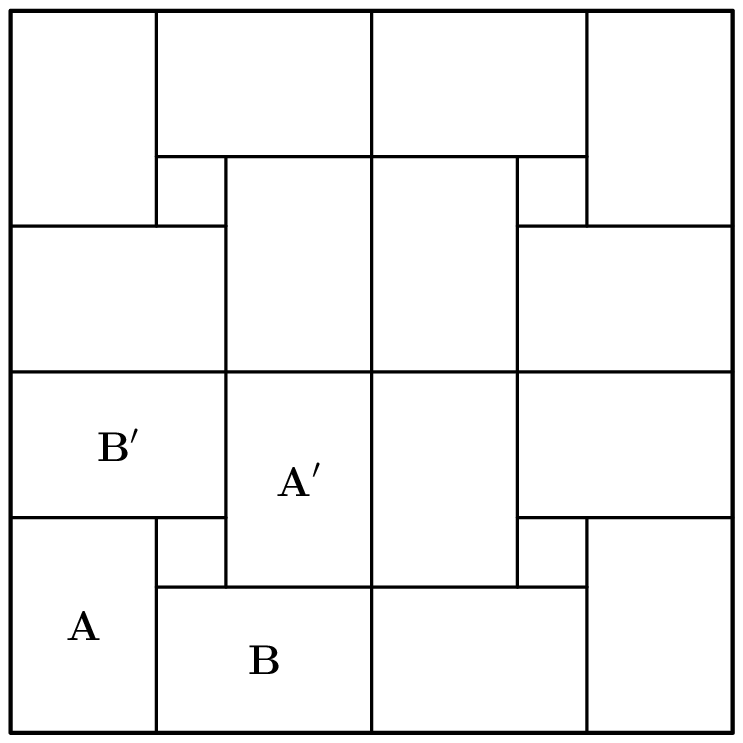} \caption{} \label{o-1} \end{figure}

\begin{figure}[!p] \centering \includegraphics{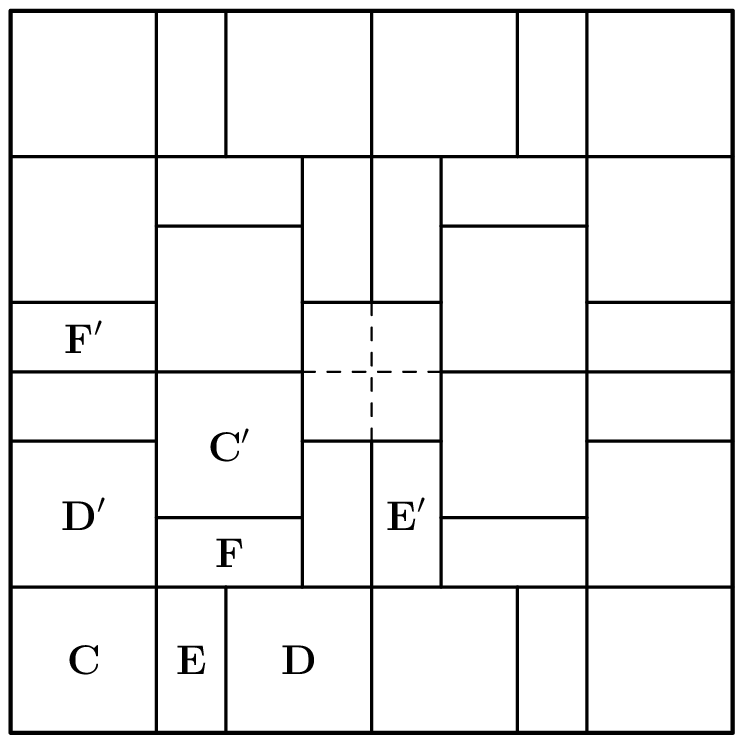} \caption{} \label{o-2} \end{figure}

\begin{figure}[!p] \centering \includegraphics{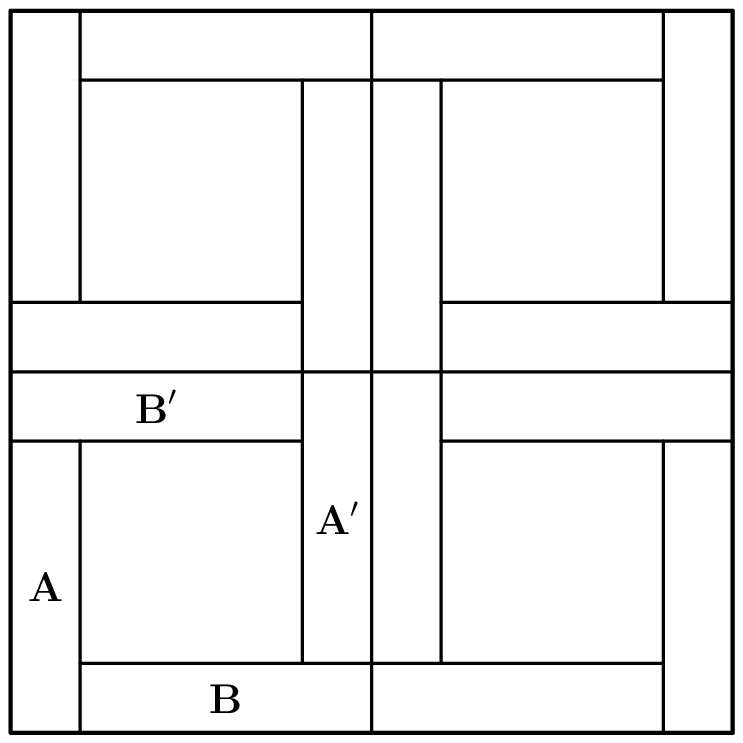} \caption{} \label{o-3} \end{figure}

\begin{figure}[!p] \centering \includegraphics{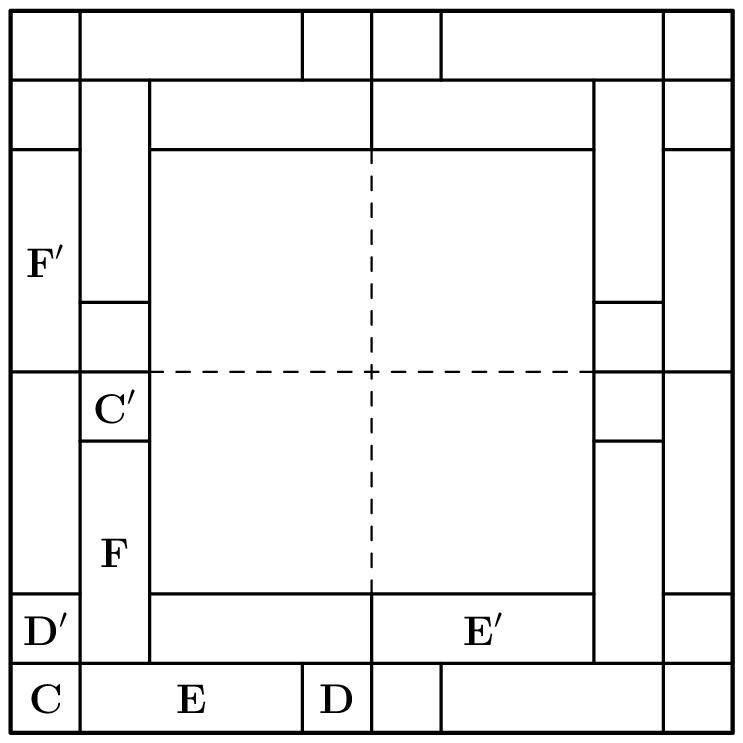} \caption{} \label{o-4} \end{figure}

\begin{figure}[!p] \centering \includegraphics{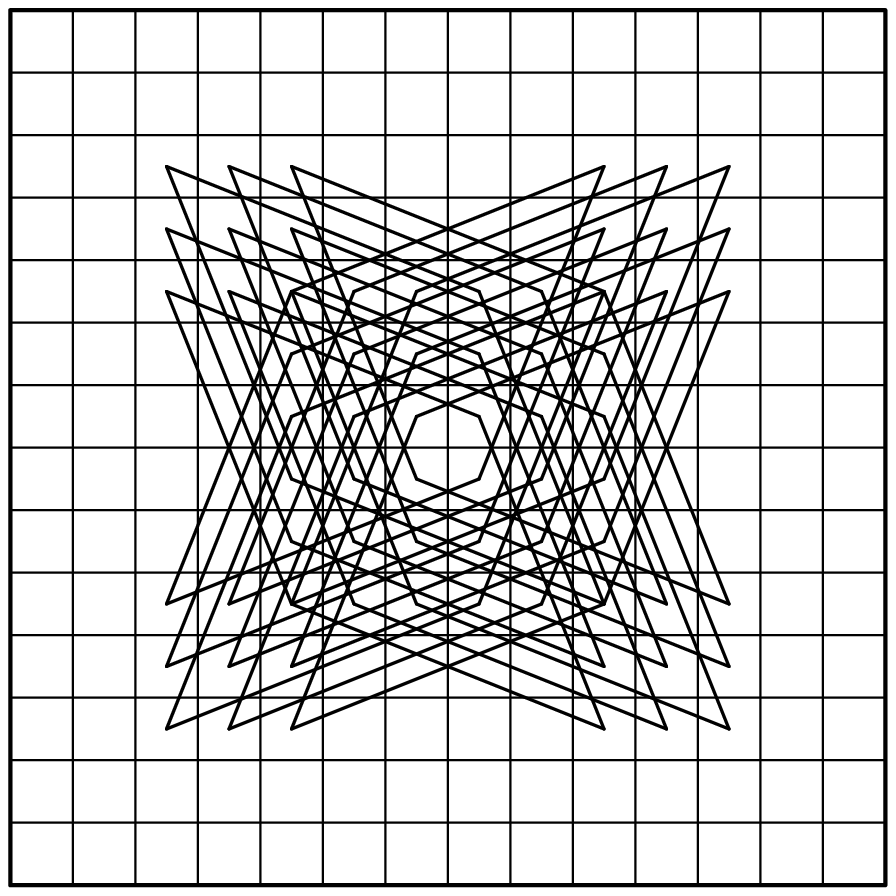} \caption{} \label{2-5-i} \end{figure}

\begin{figure}[!p] \centering \includegraphics{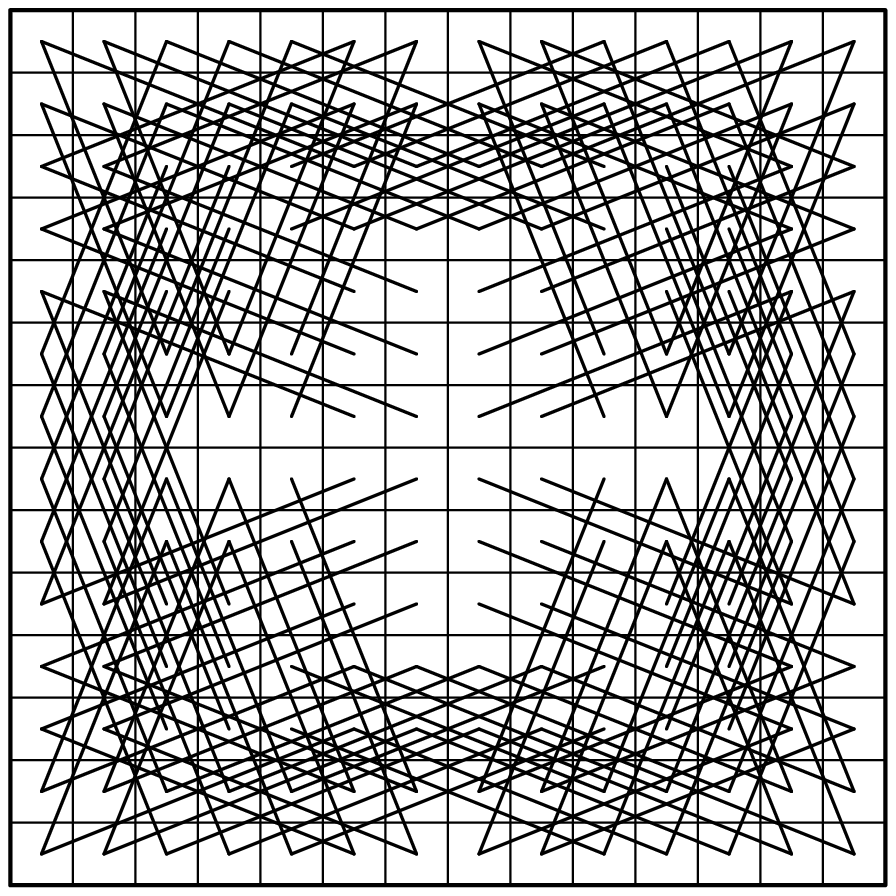} \caption{} \label{2-5-o} \end{figure}

Let the \emph{forward cores} of $B$ be the four subboards \begin{align*} C'_1 &= [p, q] \times [p, q],\\ C'_2 &= [p + q, 2q] \times [2p, p + q],\\ C'_3 &= [2p + q, p + 2q] \times [2p + q, p + 2q], \text{ and }\\ C'_4 &= [2p, p + q] \times [p + q, 2q], \end{align*} and let the \emph{backward cores} of $B$ be the four subboards \begin{align*} C''_1 &= [2p, p + q] \times [2p, p + q],\\ C''_2 &= [2p + q, p + 2q] \times [p, q],\\ C''_3 &= [p + q, 2q] \times [p + q, 2q], \text{ and }\\ C''_4 &= [p, q] \times [2p + q, p + 2q]. \end{align*}

Each core of $B$ is a square board of side $q - p$. When $2p \ge q$, all eight cores are pairwise disjoint, as in Figure \ref{i-1}. When $2p < q$, the forward core $C'_i$ and the backward core $C''_i$ overlap for all $i$, and all other pairs of cores are disjoint, as in Figure \ref{i-2}.

Consider all cycles in $W$ of length four joining corresponding cells in the four forward cores, and in the four backward cores. We refer to each such cycle as a \emph{rhombus}. The pencil of all forward rhombuses is \[C'_1 \to (q, p) \to (p, q) \to (-q, -p) \to (-p, -q),\] and the pencil of all backward rhombuses is \[C''_1 \to (q, -p) \to (-p, q) \to (-q, p) \to (p, -q).\]

The union of all rhombuses is a subgraph of $W$, the \emph{inner graph} $I$.

Furthermore, consider the six pencils \begin{align*} [0, p] \times [0, q] &\to (q, p), \tag{\textbf{A}}\\ [p, p + q] \times [0, p] &\to (-p, q), \tag{\textbf{B}}\\ [0, p] \times [0, p] &\to (p, q), \tag{\textbf{C}}\\ [q, p + q] \times [0, p] &\to (-q, p), \tag{\textbf{D}}\\ [p, q] \times [0, p] &\to (q, p), \text{ and } \tag{\textbf{E}}\\ [p, 2p] \times [p, q] &\to (-p, q) \tag{\textbf{F}} \end{align*} and their reflections. Figures \ref{o-1} and \ref{o-2} show the case $2p \ge q$, and Figures \ref{o-3} and \ref{o-4} show the case $2p < q$.

Together, all twenty-four of those pencils form a subgraph of $W$, the \emph{outer graph} $O$.

For instance, Figures \ref{2-5-i} and \ref{2-5-o} show $I$ and $O$ in the case $p = 2$ and $q = 5$.

Let the \emph{key graph} $H$ be the union of $I$ and $O$.

\medskip

\emph{Lemma 1}. The Willcocks graph admits a pseudotour.

\medskip

\emph{Proof}. Let $a$ be any cell of $B$ and let $e$ be the number of cores of $B$ that $a$ belongs to, $e = 0$, 1, or 2.

It is straightforward to verify that the degree of $a$ in $I$ is $2e$ and the degree of $a$ in $O$ is $2 - e$. Consequently, halving all rhombuses in $H$ in an arbitrary manner yields a pseudotour of $W$. $\square$

\medskip

We proceed to establish one crucial lemma.

\medskip

\emph{Lemma 2}. Let the graph $G$ be the edge-disjoint union of the cycles $C_1$, $C_2$,~\ldots, $C_k$, each of length four, and the graph $E$, in such a way that halving all of $C_1$, $C_2$,~\ldots, $C_k$ in $G$ yields a two-factor of $G$. Then $G$ is Hamiltonian if and only if it is connected.

\medskip

\emph{Proof}. Suppose that $G$ is connected.

Begin by halving all of $C_1$, $C_2$, \ldots, $C_k$ in $G$ in an arbitrary manner, obtaining a two-factor of $G$ composed of several disjoint cycles.

Consider all of $C_1$, $C_2$, \ldots, $C_k$ one by one. For each $C_i \equiv a_ib_ic_id_i$, check if its two edges in the two-factor, say $a_ib_i$ and $c_id_i$, belong to the same cycle.

If they do, leave $C_i$ as it is.

Otherwise, delete $a_ib_i$ and $c_id_i$ and replace them with $b_ic_i$ and $d_ia_i$, in effect flipping the halving of $C_i$. This operation splices the two cycles that $a_ib_i$ and $c_id_i$ belong to into a single cycle.

Once all of $C_1$, $C_2$, \ldots, $C_k$ have been considered, we obtain a two-factor $P$ of $G$ such that the two edges of each of $C_1$, $C_2$, \ldots, $C_k$ in $P$ belong to the same cycle.

Consider any edge deleted from $G$ in order to obtain $P$, say the edge $a_jb_j$ of $C_j$. Since the edges $b_jc_j$ and $d_ja_j$ belong to $P$ and they are in the same cycle, there exists a path in $P$ from $a_j$ to $b_j$.

Consequently, for every edge of $G$ there exists either an edge or a path in $P$ joining the same pair of vertices. Since $G$ is connected by supposition, so is $P$. Thus $P$ consists of a single cycle and $P$ is a Hamiltonian tour of $G$. $\square$

\medskip

Lemma 2 applies to the key graph considered as the union of all rhombuses and the outer graph. We are left to show that the key graph is in fact connected. This is the most complicated part of the proof.

To this end, first we \emph{fold} $H$ into a much simpler graph, the \emph{folding graph} $F$, of the same connectedness.

Let $q - p = 2s + 1$. The vertices of $F$ are all ordered triplets $(x, y, f)$ such that $-s \le x \le s$, $-s \le y \le s$, and $f = 1$ or $f = 2$. All vertices of the form $(x, y, 1)$ form the \emph{first floor} of $F$ and all vertices of the form $(x, y, 2)$ form the \emph{second floor} of $F$.

The edges of $F$ are defined as follows.

Given a cell $a$ that occupies position $(x + s, y + s)$ in a forward core of $B$, let $\pi_1(a)$ be the vertex $(x, y, 1)$ of $F$. Analogously, given a cell $a$ that occupies position $(x + s, y + s)$ in a backward core of $B$, let $\pi_2(a)$ be the vertex $(x, y, 2)$ of $F$.

We refer to $\pi_i(a)$ as a \emph{projection} of $a$ in $F$. When $a$ possesses a single projection (equivalently, when $a$ belongs to a single core of $B$), we denote it by $\pi(a)$.

By construction, the outer graph comprises several paths whose endpoints are precisely the cells of the symmetric difference of the forward and backward cores, possibly together with some cycles. (Eventually, we will see that the outer graph does not in fact contain any cycles.)

For each path of $O$ of endpoints $a$ and $b$, an edge in $F$ joins $\pi(a)$ and $\pi(b)$. It might happen that distinct paths in $O$ contribute the same edge to $F$; in all such cases, we consider that edge to be simple rather than doubled.

Additionally, for each cell $a$ in the intersection of two cores of $B$, an edge in $F$ joins $\pi_1(a)$ and $\pi_2(a)$. Thus, we are essentially treating all such cells as paths in $O$ of zero length.

Roughly speaking, $F$ is obtained from $H$ by replacing each path in $O$ with an edge joining its endpoints, and then collapsing each rhombus in $I$ into a single vertex so that the four forward cores are stacked on top of each other, forming the first floor of $F$, and the four backward cores are stacked on top of each other as well, forming the second floor of $F$.

\begin{figure}[!ht] \centering \includegraphics[scale=0.9]{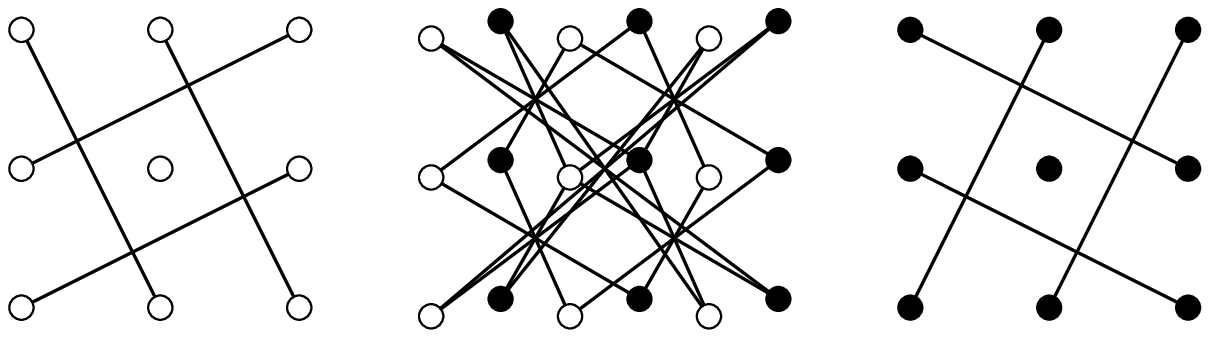} \caption{} \label{2-5-f} \end{figure}

For instance, Figure \ref{2-5-f} shows the folding graph $F$ in the case $p = 2$ and $q = 5$, corresponding to the union $H$ of the inner and outer graphs in Figures \ref{2-5-i} and \ref{2-5-o}. First-floor edges are shown on the left, second-floor edges on the right, and between-floor edges in the middle.

\medskip

\emph{Lemma 3.} If the outer graph does not contain any cycles and the folding graph is connected, then the key graph is connected as well.

\medskip

\emph{Proof}. Suppose first that the projections of two distinct cells $a'$ and $a''$ in $F$ coincide. Then $a'$ and $a''$ are connected by a path in $I$ along the edges of some rhombus.

Consider, then, any edge in $F$ joining $\pi_i(a)$ and $\pi_j(b)$. Pick $i'$, $a'$, $j'$, and $b'$ in such a way that $\pi_i(a) = \pi_{i'}(a')$, $\pi_j(b) = \pi_{j'}(b')$, and $i'$, $a'$, $j'$, and $b'$ satisfy the definition of an edge of $F$. Then there exist a path in $I$ from $a$ to $a'$, a path in $O$ (possibly of zero length) from $a'$ to $b'$, and a path in $I$ from $b'$ to $b$. Consequently, $a$ and $b$ are connected in $H$.

Lastly, let $c$ and $d$ be any two cells of $B$.

Since $O$ does not contain any cycles by supposition, there exists a path in $O$ (possibly of zero length) from $c$ to some cell $c'$ in a core of $B$. Analogously, there exists a path in $O$ from $d$ to some cell $d'$ in a core of $B$.

Since $F$ is connected by supposition, there exists a path from a projection of $c'$ to a projection of $d'$ in $F$. As above, it follows that there exists a path from $c'$ to $d'$ in $H$. $\square$

\medskip

We continue by describing one more family of graphs. Eventually, we will see that when $L$ varies over all free leapers, $F$ varies over the members of this family, and that each one of them is connected.

Let $m$ and $n$ be nonnegative integers such that $m - n$ and $m + n$ are relatively prime. The \emph{crisscross graph} $R(m, n)$ is defined as follows.

Let $m + n = 2t + 1$. The vertices of $R(m, n)$ are all ordered triplets $(x, y, f)$ such that $-t \le x \le t$, $-t \le y \le t$, and $f = 1$ or $f = 2$.

An edge joins two vertices $u_1 = (x_1, y_1, f_1)$ and $u_2 = (x_2, y_2, f_2)$ of $R(m, n)$ if and only if they satisfy one of the following conditions: \begin{align*} f_1 = f_2 = 1 &\text{ and } (x_1, y_1) \pm (m, n) = (x_2, y_2),\\ f_1 = f_2 = 1 &\text{ and } (x_1, y_1) \pm (-n, m) = (x_2, y_2),\\ f_1 = f_2 = 2 &\text{ and } (x_1, y_1) \pm (n, m) = (x_2, y_2),\\ f_1 = f_2 = 2 &\text{ and } (x_1, y_1) \pm (-m, n) = (x_2, y_2),\\ f_1 \neq f_2 &\text{ and } (x_1, y_1) \pm (\pm m, m) = (x_2, y_2), \text{ and}\\ f_1 \neq f_2 &\text{ and } (x_1, y_1) \pm (\pm n, n) = (x_2, y_2). \end{align*}

For instance, the folding graph in Figure \ref{2-5-f} coincides with $R(2, 1)$.

When $(x_1, y_1) + v = (x_2, y_2)$, we write $u_1 \to v$ for the edge in $R(m, n)$ joining $u_1$ and $u_2$, and we refer to $\pm v$ as the \emph{type} of that edge.

\medskip

\emph{Lemma 4}. The outer graph does not contain any cycles. Furthermore, let $r = q - p$, let $m$ be the common remainder of $p$ and $q$ upon division by $r$, let $n$ be determined by $m + n = r$, and let $h = \floor{\frac{p}{r}}$. Then the folding graph coincides with the crisscross graph $R(m, n)$ if $h$ is even and $R(n, m)$ if it is odd.

\medskip

\emph{Proof}. We are going to show that the outer graph is always the disjoint union of twenty pencils of paths of $L$, grouped into five families of reflections.

We distinguish three cases for $p$ and $q$.

\smallskip

\emph{Case 1}. $3p < q$. Then $m = p$, $n = q - 2p$, and $h = 0$.

\begin{figure}[!t] \centering \includegraphics{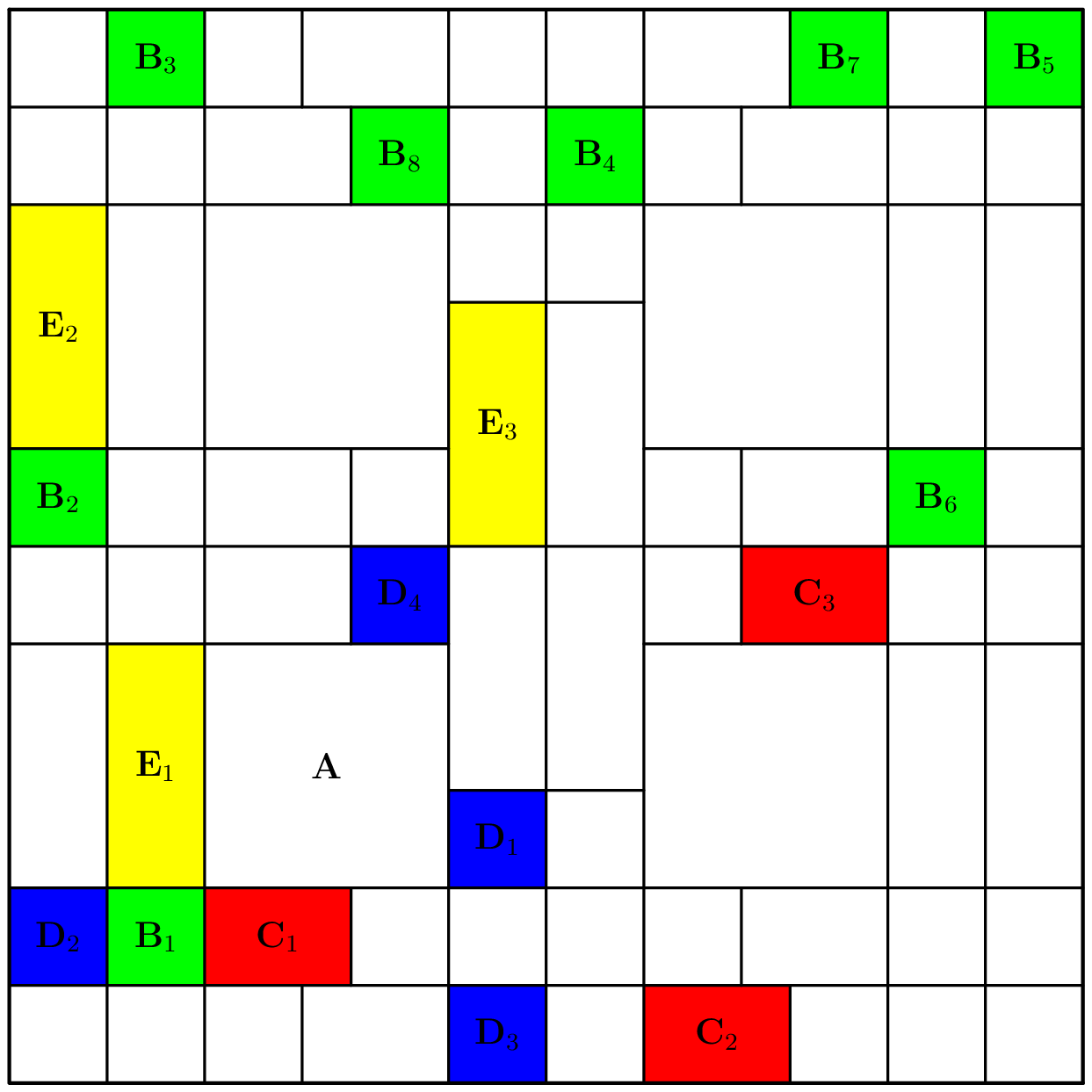} \caption{} \label{pencils-1} \end{figure}

Figure \ref{pencils-1} shows the partitioning of $O$ in this case. The twenty pencils are as follows; in all cases, we consider the given pencil together with its reflections.
~
\[[2p, q] \times [2p, q] \tag{\textbf{A}}\]

Those pencils of zero length account for all edges of types $(m, m)$ and $(-m, m)$ in $F$.
~
\begin{align*} [p, 2p] \times [p, 2p] &\to (-p, q) \to (p, q) \to (q, -p) \to (q, p) \to \tag{\textbf{B}}\\ &\to (-p, -q) \to (-p, q) \to (-q, -p) \end{align*}

Those pencils account for all edges of types $(n, n)$ and $(-n, n)$ in $F$.
~
\begin{align*} [2p, q - p] \times [p, 2p] &\to (q, -p) \to (p, q) \tag{\textbf{C}}\\ [q, p + q] \times [2p, 3p] &\to (-q, -p) \to (q, -p) \to (-p, q) \tag{\textbf{D}} \end{align*}

Those two families of pencils together account for all edges of types $(m, n)$ and $(-m, n)$ in $F$.
~
\[[p, 2p] \times [2p, q] \to (-p, q) \to (q, -p) \tag{\textbf{E}}\]

Those pencils account for all edges of types $(n, m)$ and $(n, -m)$ in $F$. 

It is straightforward to verify that all pencils we have listed do indeed fit together into $O$ as in Figure \ref{pencils-1}.

\smallskip

\emph{Case 2.} $2p \le q < 3p$. Again $m = p$, $n = q - 2p$, and $h = 0$.

\begin{figure}[!t] \centering \includegraphics{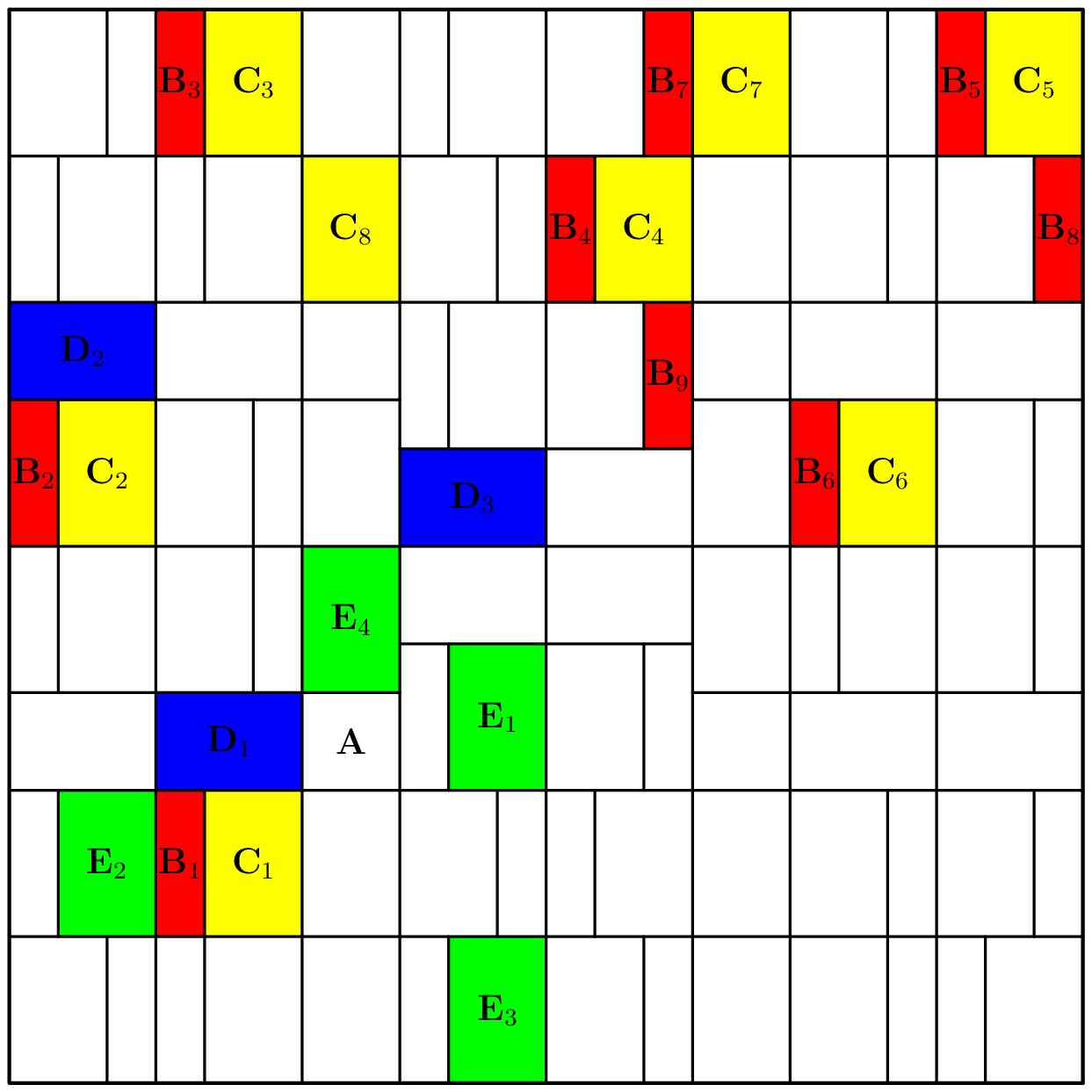} \caption{} \label{pencils-2} \end{figure}

Figure \ref{pencils-2} shows the partitioning of $O$ in this case. The twenty pencils are as follows; in all cases, we consider the given pencil together with its reflections.
~
\[[2p, q] \times [2p, q] \tag{\textbf{A}}\]

Those pencils of zero length account for all edges of types $(m, m)$ and $(-m, m)$ in $F$.
~
\begin{align*} [p, 4p - q] \times [p, 2p] &\to (-p, q) \to (p, q) \to (q, -p) \to (q, p) \to \tag{\textbf{B}}\\ &\to (-p, -q) \to (-p, q) \to (q, -p) \to (-q, -p)\\ [4p - q, 2p] \times [p, 2p] &\to (-p, q) \to (p, q) \to (q, -p) \to (q, p) \to \tag{\textbf{C}}\\ &\to (-p, -q) \to (-p, q) \to (-q, -p) \end{align*}

Those two families of pencils together account for all edges of types $(n, n)$ and $(-n, n)$ in $F$.
~
\[[p, 2p] \times [2p, q] \to (-p, q) \to (q, -p) \tag{\textbf{D}}\]

Those pencils account for all edges of types $(n, m)$ and $(-n, m)$ in $F$.
~
\[[3p, p + q] \times [2p, 3p] \to (-q, -p) \to (q, -p) \to (-p, q) \tag{\textbf{E}}\]

Those pencils account for all edges of types $(m, n)$ and $(m, -n)$ in $F$.

It is straightforward to verify that all pencils we have listed do indeed fit together into $O$ as in Figure \ref{pencils-2}.

\smallskip

\emph{Case 3.} $q < 2p$.

This case is more complicated than Cases 1 and 2 because the lengths of the pencils are no longer bounded. They cannot be, as the ratio of the number of cells in $B$ to the number of endpoints of paths of $O$ (that is, the total area of all cores of $B$) can become arbitrarily large.

For that reason, we approach it in a very different way.

We proceed by induction on $h$. Cases 1 and 2 double as the base case of the induction, $h = 0$.

Let, then, $h \ge 1$. Set $h' = h - 1$, $p' = rh' + m$, and $q' = r(h' + 1) + m$. Let $L'$ be a $(p', q')$-leaper, let $B'$ be a square chessboard of side $2(p' + q')$, and let $O'$ be the outer graph of $L'$.

Suppose that the statement of the lemma holds for $O'$. We are going to \emph{lift} $O'$ to $O$ in a way making it evident that the statement of the lemma also holds for $O$.

To this end, to each cell $a$ of $B'$ we assign a path $\varphi(a)$ in $O$, possibly of zero length, so that the following conditions hold.

(a) The cells of $B$ are the disjoint union of the vertex sets of the paths $\varphi(a)$ for $a$ in $B'$.

(b) If $a$ occupies position $(x, y)$ in a forward or backward core of $B'$, then an endpoint of $\varphi(a)$ occupies position $(x, y)$ in a backward or forward core of $B$, respectively.

(c) If an edge of $O'$ joins $a$ and $b$, then an edge of $O$ joins an endpoint of $\varphi(a)$ and an endpoint of $\varphi(b)$.

Since every vertex of $O$ in a core of $B$ is of degree one and every other vertex of $O$ is of degree two, the two endpoints of condition (c) are always uniquely determined. Consequently, conditions (a), (b), and (c) ensure that all paths $\varphi(a)$ for $a$ in $B'$ join together to form all of $O$.

Since by the induction hypothesis $O'$ does not contain any cycles, neither does $O$.

Since toggling the floor of every vertex transforms each of the crisscross graphs $R(m, n)$ and $R(n, m)$ into the other, by the induction hypothesis and condition (b) it follows that $F$ does indeed coincide with a crisscross graph as required and the second part of the lemma holds for $O$ as well.

The mapping $\varphi$ takes on two different forms in the cases $h = 1$ and $h \ge 2$.

\begin{figure}[!p] \centering \includegraphics{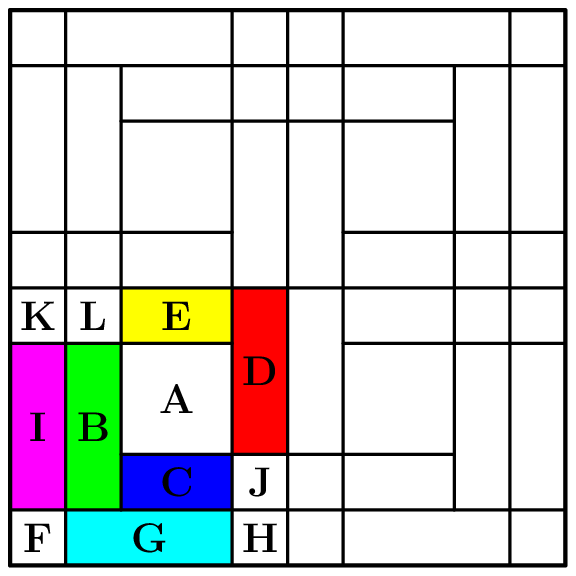} \caption{} \label{lift-1a} \end{figure}

\begin{figure}[!p] \centering \includegraphics{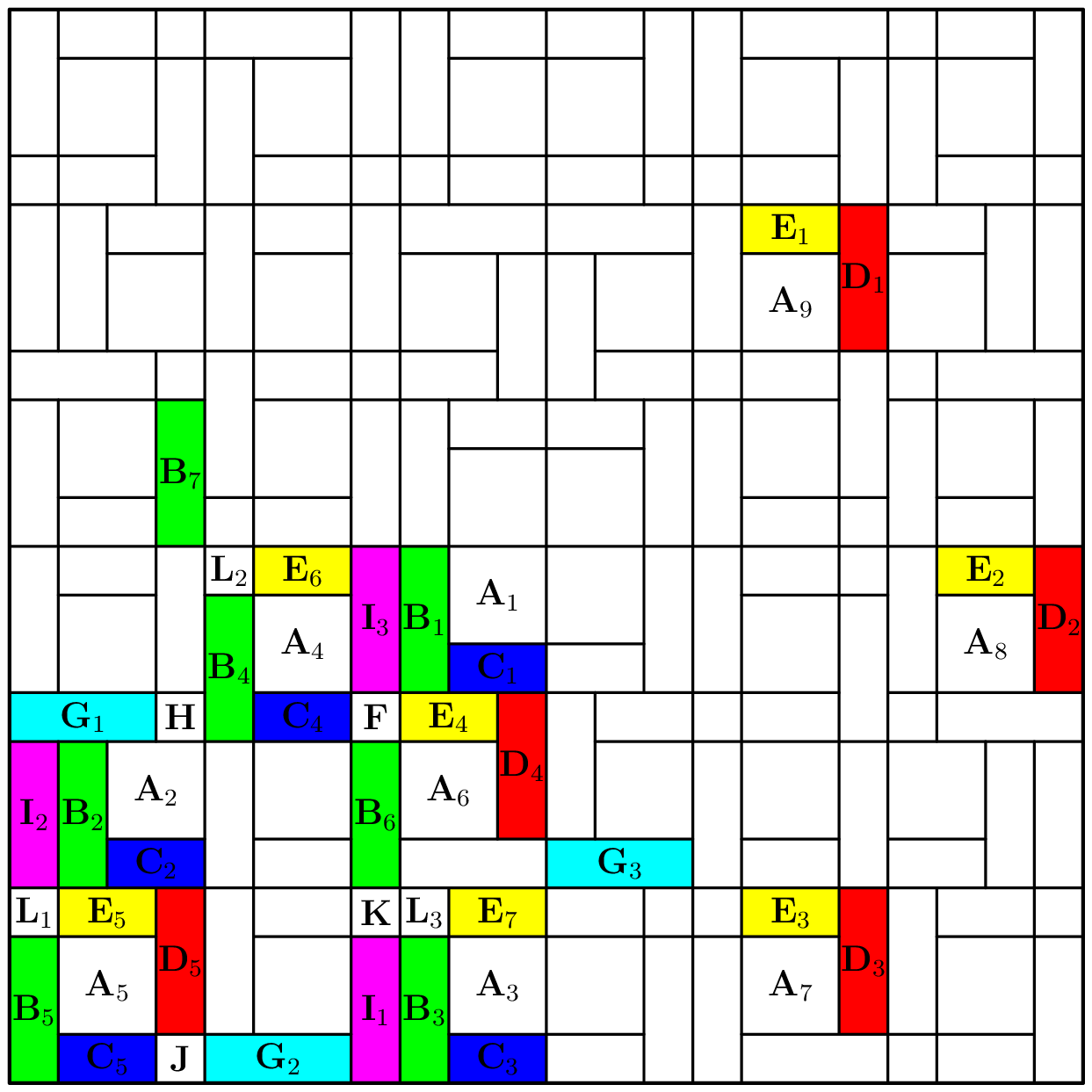} \caption{} \label{lift-1b} \end{figure}

When $h = 1$, $\varphi$ is as follows.

The expression $S': S \to d_1 \to d_2 \to \ldots \to d_k$ indicates that if $a$ occupies position $(x, y)$ in a subboard $S'$ of $B'$ and $b$ is the cell at position $(x, y)$ in the subboard $S$ of $B$, then $\varphi(a) = b \to d_1 \to d_2 \to \ldots \to d_k$. By symmetry, the definition of $\varphi$ extends to all reflections of the subboards $S'$ in $B'$.
~
\begin{align*}
% A
[2p', q'] \times [2p', q'] : {}&[p' + 2q', 3q' - p'] \times [p' + 2q', 3q' - p'] \to \tag{\textbf{A}}\\ &\to (-q, -p) \to (q, -p) \to (-p, q) \to (-p, -q) \to\\ &\to (q, p) \to (q, -p) \to (p, q) \to (-p, q)\\
% B
[p', 2p'] \times [p', q'] : {}&[2q', p' + 2q'] \times [2q', 3q' - p'] \to \tag{\textbf{B}}\\ &\to (-q, -p) \to (q, -p) \to (-p, q) \to (-p, -q) \to\\ &\to (q, p) \to (-p, q)\\
% C
[2p', q'] \times [p', 2p'] : {}&[p' + 2q', 3q' - p'] \times [2q', 2q' + p'] \to \tag{\textbf{C}}\\ &\to (-q, -p) \to (q, -p) \to (-p, q) \to (-p, -q)\\
% D
[q', p' + q'] \times [2p', p' + q'] : {}&[5q' - 3p', 5q' - 2p'] \times [4q' - p', 5q' - 2p'] \to \tag{\textbf{D}}\\ &\to (p, -q) \to (-p, -q) \to (-q, p) \to (-q, -p)\\
% E
[2p', q'] \times [q', p' + q'] : {}&[4q' - p', 5q' - 3p'] \times [5q' - 3p', 5q' - 2p'] \to \tag{\textbf{E}}\\ &\to (p, -q) \to (-p, -q) \to (-q, p) \to (-q, -p) \to\\ &\to (p, q) \to (p, -q)\\
% F
[0, p'] \times [0, p'] : {}&[2q' - p', 2q'] \times [2q' - p', 2q'] \tag{\textbf{F}}\\
% G
[p', q'] \times [0, p'] : {}&[0, q' - p'] \times [2q' - p', 2q'] \to \tag{\textbf{G}}\\ &\to (p, -q) \to (q, p)\\
% H
[q', p' + q'] \times [0, p'] : {}&[q' - p', q'] \times [2q' - p', 2q'] \tag{\textbf{H}}\\
% I
[0, p'] \times [p', q'] : {}&[2q' - p', 2q'] \times [0, q' - p'] \to \tag{\textbf{I}}\\ &\to (-q, p) \to (q, p)\\
% J
[q', p' + q'] \times [p', 2p'] : {}&[q' - p', q'] \times [0, p'] \tag{\textbf{J}}\\
% K
[0, p'] \times [q', p' + q'] : {}&[2q' - p', 2q'] \times [q' - p', q'] \tag{\textbf{K}}\\
% L
[p', 2p'] \times [q', p' + q'] : {}&[0, p'] \times [q' - p', q'] \to \tag{\textbf{L}}\\ &\to (p, q) \to (p, -q)
~
\end{align*}

Figure \ref{lift-1a} shows the partitioning of $B'$ into subboards, and Figure \ref{lift-1b} shows the pencils that the paths assigned by $\varphi$ to the cells of each subboard form in $B$.

(Alternatively, instead of giving this form of $\varphi$ we could have explicitly described all twenty pencils in the cases $h = 1$ and $m < n$, and $h = 1$ and $m > n$, much as we did in Cases 1 and 2.)

\begin{figure}[!p] \centering \includegraphics{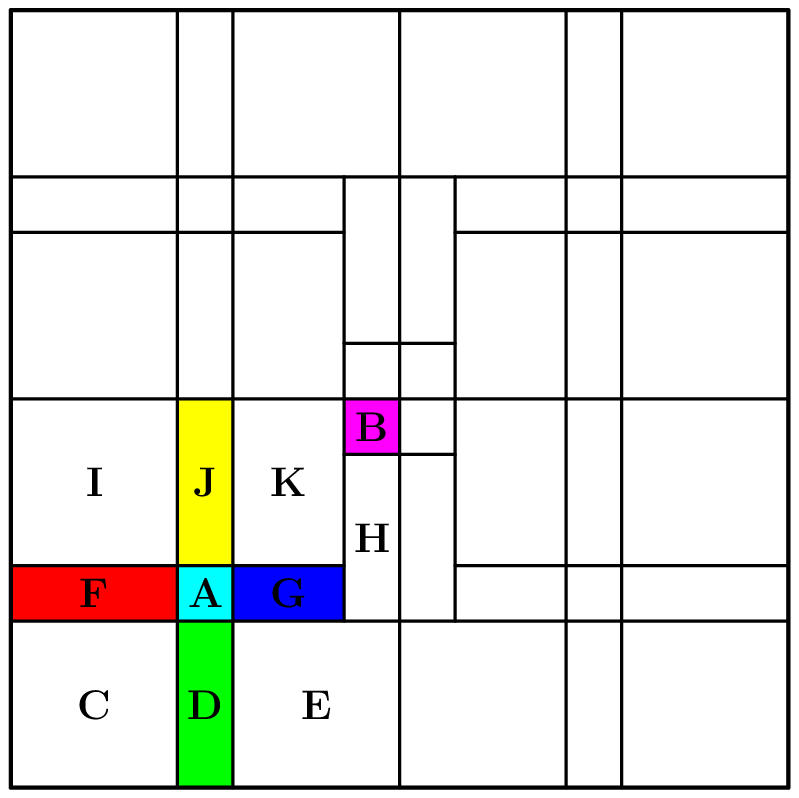} \caption{} \label{lift-2a} \end{figure}

\begin{figure}[!p] \centering \includegraphics{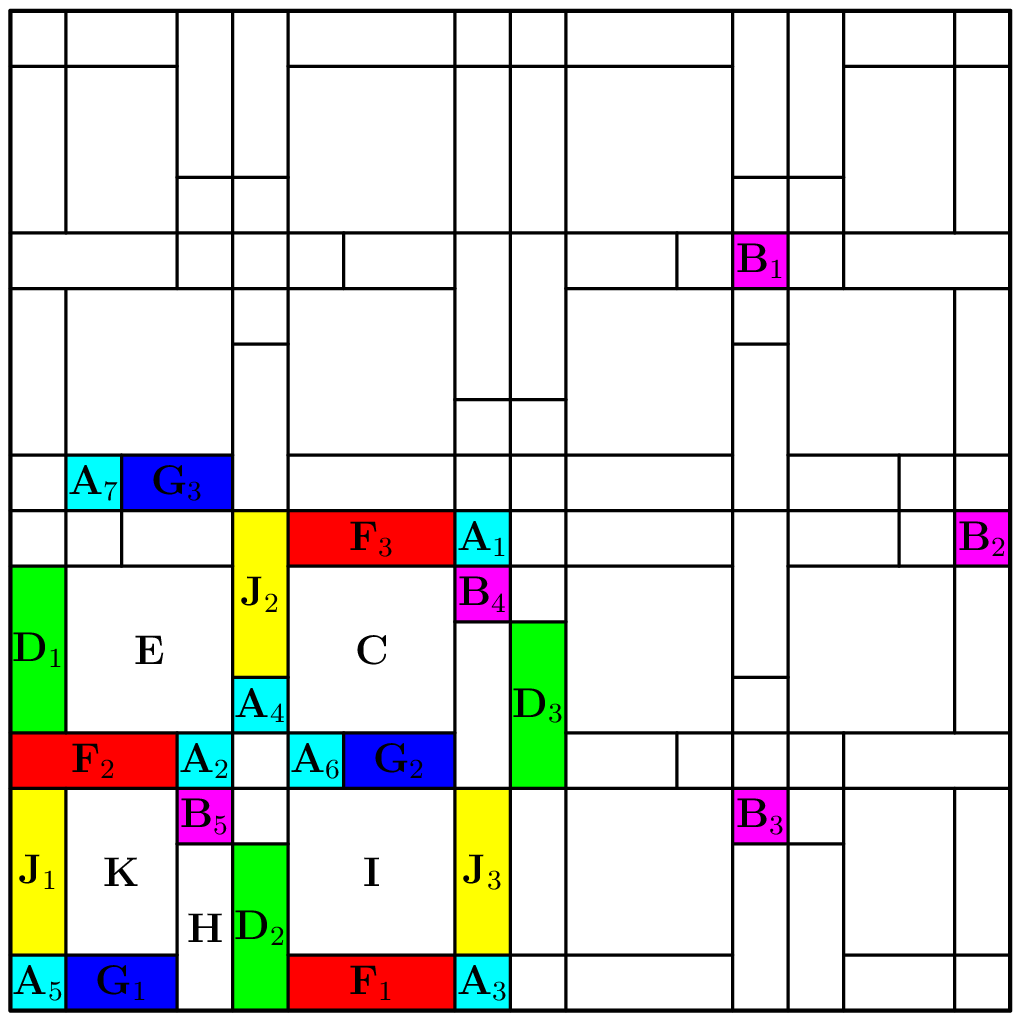} \caption{} \label{lift-2b} \end{figure}

When $h \ge 2$, the mapping $\varphi$ is as follows.
~
\begin{align*}
% A
[p', q'] \times [p', q'] : {}&[2q', 3q' - p'] \times [2q', 3q' - p'] \to \tag{\textbf{A}}\\ &\to (-q, -p) \to (q, -p) \to (-p, q) \to (-p, -q) \to\\ &\to (q, p) \to (-p, q)\\
% B
[2p', p' + q'] \times [2p', p' + q'] : {}&[4q' - p', 5q' - 2p'] \times [4q' - p', 5q' - 2p'] \to \tag{\textbf{B}}\\ &\to (p, -q) \to (-p, -q) \to (-q, p) \to (-q, -p)\\
% C
[0, p'] \times [0, p'] : {}&[2q' - p', 2q'] \times [2q' - p', 2q'] \tag{\textbf{C}}\\
% D
[p', q'] \times [0, p'] : {}&[0, q' - p'] \times [2q' - p', 2q'] \to \tag{\textbf{D}}\\ &\to (p, -q) \to (q, p)\\
% E
[q', p' + q'] \times [0, p'] : {}&[q' - p', q'] \times [2q' - p', 2q'] \tag{\textbf{E}}\\
% F
[0, p'] \times [p', q'] : {}&[2q' - p', 2q'] \times [0, q' - p'] \to \tag{\textbf{F}}\\ &\to (-q, p) \to (q, p)\\
% G
[q', 2p'] \times [p', q'] : {}&[q' - p', p'] \times [0, q' - p'] \to \tag{\textbf{G}}\\ &\to (q, p) \to (-p, q)\\
% H
[2p', p' + q'] \times [p', 2p'] : {}&[p', q'] \times [0, p'] \tag{\textbf{H}}\\
% I
[0, p'] \times [q', p' + q'] : {}&[2q' - p', 2q'] \times [q' - p', q'] \tag{\textbf{I}}\\
% J
[p', q'] \times [q', p' + q'] : {}&[0, q' - p'] \times [q' - p', q'] \to \tag{\textbf{J}}\\ &\to (p, q) \to (p, -q)\\
% K
[q', 2p'] \times [q', p' + q'] : {}&[q' - p', p'] \times [q' - p', q'] \tag{\textbf{K}}
~
\end{align*}

Figure \ref{lift-2a} shows the partitioning of $B'$ into subboards, and Figure \ref{lift-2b} shows the pencils that the paths assigned by $\varphi$ to the cells of each subboard form in $B$.

The verification of all three conditions (a), (b), and (c) for both forms of $\varphi$ is straightforward. $\square$

\medskip

\emph{Lemma 5}. Every crisscross graph is connected.

\medskip

\emph{Proof.} Since, as above, toggling the floor of every vertex switches around the crisscross graphs $R(m, n)$ and $R(n, m)$, it suffices to consider the case $m < n$.

We proceed by induction on $m + n$.

When $m = 0$ and $n = 1$, $R(0, 1)$ consists of two vertices joined by an edge and is thus connected.

Let $0 < m < n$. We distinguish three cases for $m$ and $n$.

\smallskip

\emph{Case 1}. $3m < n$. Set $m' = m$ and $n' = n - 2m$. Then $m' < n'$, $m' - n'$ and $m' + n'$ are relatively prime, and $m' + n' < m + n$.

Notice that the vertex set of $R(m', n')$ is a subset of the vertex set of $R(m, n)$.

First we show that, for every vertex $u = (x, y, f)$ of $R(m, n)$, there exists a path in $R(m, n)$ from $u$ to a vertex of $R(m', n')$ of length at most two.

If $-\floor{\frac{1}{2}(n - m)} \le x \le \floor{\frac{1}{2}(n - m)}$ and $-\floor{\frac{1}{2}(n - m)} \le y \le \floor{\frac{1}{2}(n - m)}$, then $u$ is already a vertex of $R(m', n')$.

If $\ceil{\frac{1}{2}(n - m)} \le x$ and $\ceil{\frac{1}{2}(n - m)} \le y$, then the edge $u \to (-m, -m)$ of $R(m, n)$ connects $u$ to a vertex of $R(m', n')$.

If $-\floor{\frac{1}{2}(n - m)} \le x \le \floor{\frac{1}{2}(n - m)}$, $\ceil{\frac{1}{2}(n - m)} \le y$, and $f = 1$, then the path $u \to (-m, -n) \to (m, m)$ in $R(m, n)$ leads from $u$ to a vertex of $R(m', n')$.

All other cases for $u$ are handled symmetrically.

Then we show that, for every edge $e$ of $R(m', n')$ joining $u_1 = (x_1, y_1, f_1)$ and $u_2 = (x_2, y_2, f_2)$, there exists a path in $R(m, n)$ connecting $u_1$ and $u_2$ of length either one or three.

If $e = u_1 \to (m', m')$, then $e$ is already an edge of $R(m, n)$.

If $e = u_1 \to (n', n')$, then the path $u_1 \to (-m, -m) \to (n, n) \to (-m, -m)$ in $R(m, n)$ leads from $u_1$ to $u_2$.

If $f_1 = f_2 = 1$ and $e = u_1 \to (m', n')$, then the path $u_1 \to (m, -m) \to (-m, n) \to (m, -m)$ in $R(m, n)$ leads from $u_1$ to $u_2$.

All other cases for $e$ are handled symmetrically.

Let, then, $w_1$ and $w_2$ be any two vertices of $R(m, n)$.

By our first observation, there exist a path in $R(m, n)$ from $w_1$ to some vertex $w'_1$ of $R(m', n')$ and a path in $R(m, n)$ from $w_2$ to some vertex $w'_2$ of $R(m', n')$.

By the induction hypothesis, $R(m', n')$ is connected. Consequently, there exists a path in $R(m', n')$ leading from $w'_1$ to $w'_2$. By our second observation, it follows that there exists a path from $w'_1$ to $w'_2$ in $R(m, n)$ as well.

Therefore, the crisscross graph $R(m, n)$ is connected.

\smallskip

\emph{Case 2}. $2m \le n < 3m$. Set $m' = n - 2m$ and $n' = m$. Then $m' < n'$, $m' - n'$ and $m' + n'$ are relatively prime, and $m' + n' < m + n$.

From this point on, Case 2 is fully analogous to Case 1.

\smallskip

\emph{Case 3}. $n < 2m$. Set $m' = 2m - n$ and $n' = m$. Then $m' < n'$, $m' - n'$ and $m' + n'$ are relatively prime, and $m' + n' < m + n$.

From this point on, Case 3 is fully analogous to Cases 1 and 2. $\square$

\medskip

By Lemmas 1--5, we obtain the following theorem.

\medskip

\emph{Theorem 1}. (Willcocks's Conjecture) The Willcocks graph admits a Hamiltonian tour.

\section{Symmetric Tours and Larger Boards}

The Hamiltonian tours given by Theorem 1 are not, in general, symmetric. We proceed to show how to augment the algorithm of Lemma 2 so that it yields symmetric tours.

\medskip

\emph{Theorem 2}. The Willcocks graph admits a centrally symmetric Hamiltonian tour.

\medskip

\emph{Proof.} We are going to use ``centrally symmetric'' as shorthand for ``symmetric with respect to the center of $B$''.

Start by halving all rhombuses in a centrally symmetric manner. This turns $H$ into a centrally symmetric pseudotour $P$ of $W$ composed of several disjoint cycles.

Consider the unique forward rhombus $r_1 = a_1b_1c_1d_1$ symmetric with respect to the center of $B$. If the two edges of $r_1$ in $P$ belong to different cycles, flip the halving of $r_1$ as in the proof of Lemma 2. This ensures that all edges of $r_1$ in $P$ belong to the same centrally symmetric cycle $C$.

Let $s'$ and $s''$ be the two paths obtained from $C$ by removing the edges of $r_1$ in $P$. The path $s'$ cannot be centrally symmetric, as it would then have to contain either a cell or an edge symmetric with respect to the center of $B$, and none such exist. Therefore, $s'$ and $s''$ are reflections of each other in the center of $B$. Consequently, the cells $a_1$, $b_1$, $c_1$, and $d_1$ occur along $C$ in this order.

Suppose that $C$ contains one of the edges of some rhombus $a'b'c'd'$ in $P$, say $a'b'$, but not its other edge in $P$, $c'd'$.

Let the rhombus $a''b''c''d''$ be the reflection of $a'b'c'd'$ in the center of $B$. Then the edge $a''b''$ of $a''b''c''d''$ belongs to $C$ and its edge $c''d''$ belongs to $P$ but not to $C$. Moreover, all eight vertices of $a'b'c'd'$ and $a''b''c''d''$ are pairwise distinct.

Without loss of generality, the cells $a_1$, $a'$, $b'$, $b_1$, $c_1$, $a''$, $b''$, and $d_1$ occur along $C$ in this order, possibly with $a_1$ and $a'$, $b'$ and $b_1$, $c_1$ and $a''$, or $b''$ and $d_1$ coinciding.

Suppose first that the edges $c'd'$ and $c''d''$ belong to different cycles in $P$. Flip the halvings of $a'b'c'd'$ and $a''b''c''d''$. This operation splices the cycles of $r_1$, $c'd'$, and $c''d''$ into a single centrally symmetric cycle containing all edges of $a'b'c'd'$ and $a''b''c''d''$ in the altered pseudotour.

Suppose, then, that $c'd'$ and $c''d''$ belong to the same cycle $D$ in $P$. As above, the cells $c'$, $d'$, $c''$, and $d''$ occur along $D$ in this order. Flip the halvings of $r_1$, $a'b'c'd'$, and $a''b''c''d''$. Again, this operation splices the cycle of $r_1$ and $D$ into a single centrally symmetric cycle containing all edges of $a'b'c'd'$ and $a''b''c''d''$ in the altered pseudotour.

\begin{figure}[!t] \centering \includegraphics{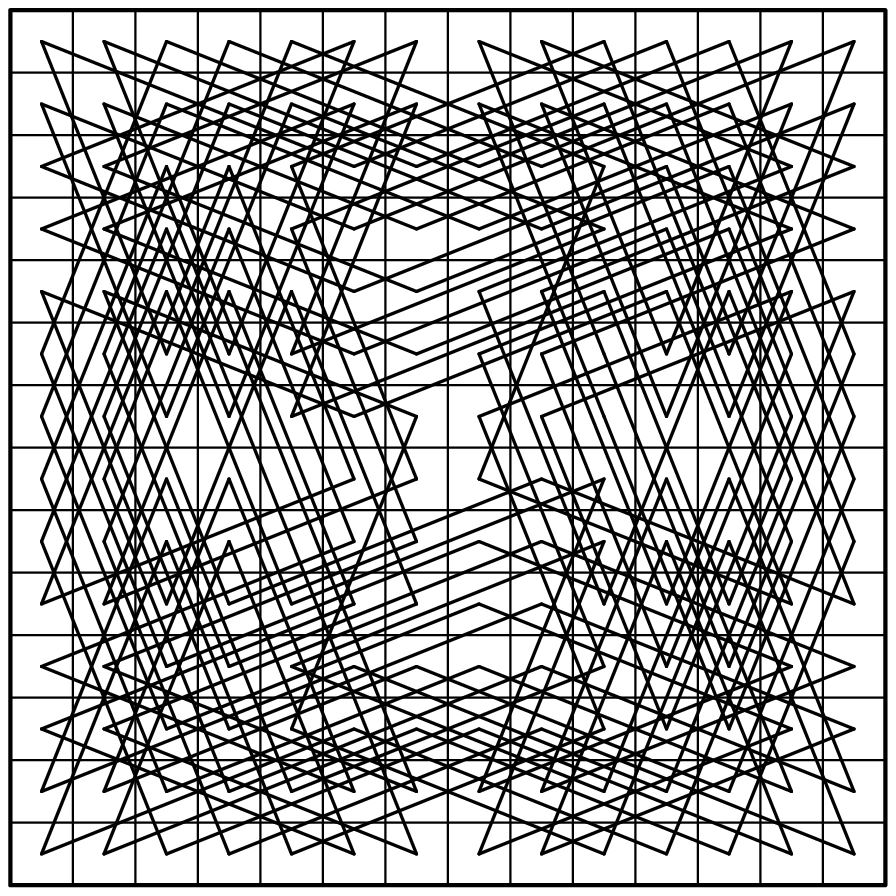} \caption{} \label{2-5-hs} \end{figure}

\begin{figure}[!p] \centering \includegraphics{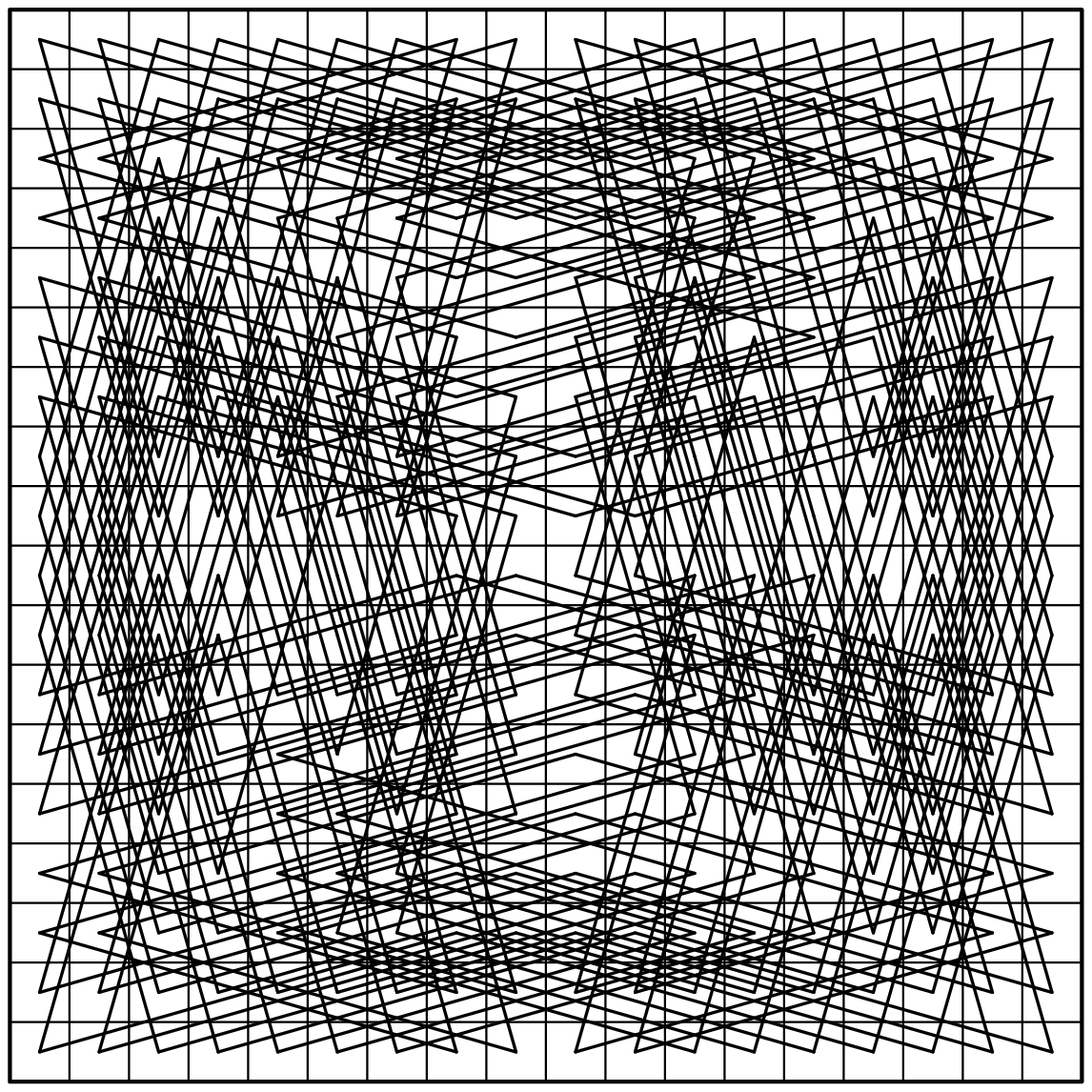} \caption{} \label{2-7-hs} \end{figure}

\begin{figure}[!p] \centering \includegraphics{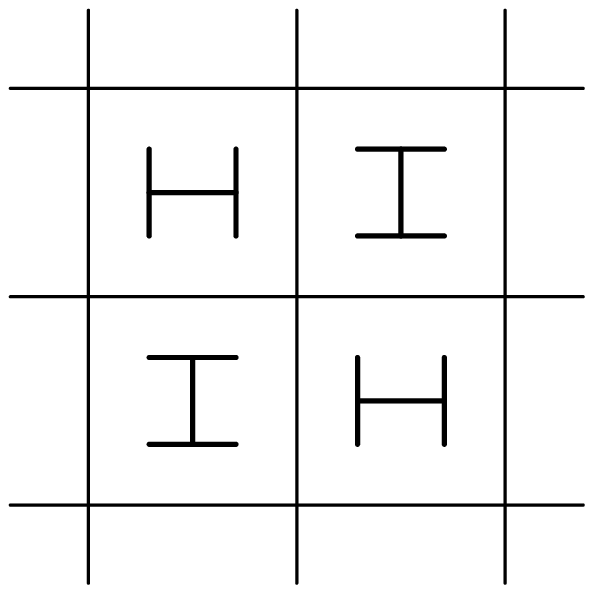} \caption{} \label{large} \end{figure}

Continue to expand $C$ in this way until for every rhombus such that one of its edges in $P$ is in $C$, its other edge in $P$ is also in $C$. We are going to show that at this point $C$ contains all cells of $B$ and is thus a centrally symmetric Hamiltonian tour of $W$.

To this end, let $a$ be any cell of $B$. Since $H$ is connected, there exists a path from a vertex of $r_1$ to $a$ in $H$. Let $r_2$, $r_3$, \ldots, $r_k$ be rhombuses and $p_1$, $p_2$, \ldots, $p_k$ paths in the outer graph $O$ (some possibly of zero length) such that, for all $i < k$, $p_i$ connects a vertex of $r_i$ to a vertex of $r_{i + 1}$, and $p_k$ connects a vertex of $r_k$ to $a$.

Since $C$ contains all vertices of $r_1$, it contains the path $p_1$. Thus $C$ contains at least one vertex of $r_2$, namely the endpoint of $p_1$ in $r_2$. By the special property of $C$, since it contains at least one vertex of $r_2$ it must contain all of them. Therefore, $C$ contains the path $p_2$.

By going on in this way, eventually we obtain that $a$ belongs to $C$. $\square$

\medskip

Figures \ref{2-5-hs} and \ref{2-7-hs} show centrally symmetric Hamiltonian tours of the $(2, 5)$-leaper and the $(2, 7)$-leaper obtained by the algorithm of Theorem 2.

We conclude by extending Theorems 1 and 2 to arbitrarily large boards.

\medskip

\emph{Theorem 3}. There exists a Hamiltonian tour of $L$ on every rectangular board both of whose sides are multiples of $2(p + q)$.

\medskip

\emph{Proof}. Dissect the board into square subboards of side $2(p + q)$.

Let $C$ be any Hamiltonian tour of $B$ given by Theorems 1 or 2. Within each square subboard, place either a translation copy of $C$ or a copy of $C$ rotated by $90^\circ$, in such a way that copies adjacent by side are always of different orientations as in Figure \ref{large}.

We say that two edges $ab$ and $cd$ in neighbouring copies of $C$ form a \emph{switch} if both of $bc$ and $da$ are moves of $L$. Flipping a switch, as in the proofs of Lemma 2 and Theorem 2, splices the cycles of $ab$ and $cd$ into a single cycle.

Consider two neighbouring copies of $C$, a translation one $C'$ in the subboard $B'$ on the left and a rotation one $C''$ in the subboard $B''$ on the right.

Let $a$ and $b$ be the cells at positions $(2p + q, 0)$ and $(3p + q, q)$ in $B'$. Since $ab$ is an edge in the outer graph of $B'$, it is also an edge in $C'$. Analogously, the cells $c$ and $d$ at positions $(p, p + q)$ and $(0, p)$ in $B''$ are the endpoints of an edge in $C''$. Thus $ab$ and $cd$ form a switch.

By symmetry, there exists a switch between each pair of neighbouring copies of $C$. Moreover, no two of those switches share an edge.

Consequently, flipping several switches allows us to splice all copies of $C$ into a Hamiltonian tour of the board. $\square$

\medskip

When at least one side of the board is the product of $2(p + q)$ and an odd positive integer, it is possible to arrange the flipped switches in a way symmetric with respect to the center of the board. By further picking $C$ to be centrally symmetric as well, we obtain the following corollary.

\medskip

\emph{Corollary}. There exist centrally symmetric Hamiltonian tours of $L$ on arbitrarily large square boards.

\section{Acknowledgements}

I am grateful to Professor Donald E. Knuth for reading a preprint of this paper and contributing extensive comments and valuable advice.

\end{document}